\documentclass[11pt,a4paper,leqno]{article}
\usepackage{amsmath}
\usepackage{amssymb}
\usepackage{graphicx}
\advance\topmargin by -2.2cm
\newcommand{\bbt}{T\!\!\!T}
\newcommand{\bbr}{I\!\!R}

\newcommand{\bbn}{I\!\!N}

\newcommand{\calc}{{\cal C}}

\newcommand{\call}{{\cal L}}

\newcommand{\calo}{{\cal O}}

\newcommand{\barr}{\begin{array}}
\newcommand{\earr}{\end{array}}
\newcommand{\beqq}{\begin{equation}}
\newcommand{\eeqq}{\end{equation}}
\newcommand{\beao}{\begin{eqnarray*}}
\newcommand{\eeao}{\end{eqnarray*}\noindent}
\newcommand{\beam}{\begin{eqnarray}}
\newcommand{\eeam}{\end{eqnarray}\noindent}

\newcommand{\halmos}{\quad\hfill\mbox{$\Box$}}

\newcommand{\la}{\lambda}

\newcommand{\al}{\alpha}

\newcommand{\om}{\omega}

\newcommand{\vep}{\varepsilon}

\newtheorem{theo}{Theorem}
\newtheorem{prop}{\indent Proposition}

\newtheorem{rem}{\indent Remark}

\newtheorem{defin}{\indent Definition}
\newtheorem{cor}{\indent Corollary}

\newtheorem{ex}{\indent Example}

\newcommand{\wt}{\widetilde}

\newcommand{\ov}{\overline}

\newcommand{\Lra}{\Longrightarrow}

\newcommand{\lra}{\longrightarrow}

\newcommand{\nto}{n\to\infty}

\author{R. H\"opfner$^*$ \and E.~L\"ocherbach \and M. Thieullen\thanks{This work has been supported by the Agence Nationale de la Recherche through the project MANDy, Mathematical Analysis of Neuronal Dynamics, ANR-09-BLAN-0008-01. e-mail addresses: \tt{hoepfner@mathematik.uni-mainz.de},
     \tt{eva.loecherbach@u-cergy.fr} and
    \tt{michele.thieullen@upmc.fr}} \\
{\it Johannes Gutenberg-Universit\"at Mainz, Universit\'e de Cergy-Pontoise} \\
 {\it and Universit\'e Pierre et Marie Curie.}}
\begin{document}

\title{Ergodicity for a stochastic Hodgkin-Huxley model driven by Ornstein-Uhlenbeck type input}

\maketitle

\begin{abstract}
We consider a model describing a neuron and the input it receives from its dendritic tree when this input is a random perturbation of a periodic deterministic signal, driven by an Ornstein-Uhlenbeck process. The neuron itself is modeled by a variant of the classical Hodgkin-Huxley model. Using the existence of an accessible point where the weak H\"ormander condition holds and the fact that the coefficients of the system are analytic, we show that the system is non-degenerate. The existence of a Lyapunov function allows to deduce the existence of (at most a finite number of) extremal invariant measures for the process. As a consequence, the complexity of the system is drastically reduced in comparison with the deterministic system.
\end{abstract}

\\


\section{Introduction}
In this paper we study a stochastic model for a spiking neuron together with the input it receives from its dendritic tree. Our model is derived from the well-known deterministic Hodgkin-Huxley model and takes the form of a highly degenerate time inhomogeneous stochastic system.

The deterministic Hodgkin-Huxley model for the membrane potential of a neuron has been extensively studied over the last decades. 
There seems to be a large agreement (see e.g.\ the introduction in Destexhe 1997) that the $4$-dimensional dynamical system proposed initially by Hodgkin and Huxley 1952 models adequately the mechanism of spike generation in response to an external stimulus in many types of neurons. Hodgkin and Huxley modeled the behavior of ion channels with respect to the two ion currents which are predominant (import of Na$^+$ and export of K$^+$ ions through the membrane) in a way which later was found experimentally (cf.\  Izhikevich 2007, Figure 2.8 on p.\ 33) to correspond to a structure in the voltage gated ion channels which was not yet observable in 1952. Also generalizations of deterministic Hodgkin-Huxley models taking into account a larger number of types of ion channels have been considered; for a modern introduction see Izhikevich 2007.

A classical deterministic Hodgkin-Huxley system (see e.g.\ Izhikevich 2007, pp.\ 33 and 37--38) has four variables, the voltage (measured by some electrode in the soma of the neuron) and three gating variables (the state of specific voltage sensors which activate or deactivate ion channels). In addition there is some fixed deterministic function of time which represents an input. Our stochastic Hodgkin-Huxley model has only one source of stochasticity: we are interested in the effect of an external noise on the behavior of the system. Thus we replace deterministic input by the increments of a stochastic process whose stochastic differential equation plays the role of a fifth equation. A cortical neuron belonging to an active cortical network receives its input from a large number of other neurons through a huge number of synapses located on a dendritic tree of complex topological structure. This is our reason for modeling dendritic input as an autonomous diffusion process $(\xi_t)_{t\ge 0}$, time inhomogeneous and of mean-reverting Ornstein-Uhlenbeck type, having some $T$-periodic deterministic signal $t\to S(t)$ coded in its semigroup. We think of $t\to S(t)$ as a signal processed by the network: roughly speaking the signal is present in the mean values of the diffusion process as a function of time $t\ge 0$. For the three gating variables we keep the corresponding Hodgkin-Huxley equations unchanged: their activity is conditionally deterministic given the voltage, without intrinsic source of randomness. Our equation for the voltage is keeping the traditional Hodgkin-Huxley form of the drift coefficient (a function of the voltage and the gating variables) but replaces the classical deterministic input by a stochastic input $\,d\xi_t\,$ at time $t\ge 0$. In this way, we are led to consider a $5$-dimensional random dynamical system $(X_t)_{t\ge 0}$ governed by one-dimensional Brownian motion which represents the external noise: the driving Brownian motion of the Ornstein-Uhlenbeck type SDE, present in two of its five equations. 

The present paper is the second part of our study of periodic ergodicity for such models. The first part is the companion paper H\"opfner, L\"ocherbach and Thieullen 2013 where we address the existence of densities for strongly degenerate time inhomogeneous random models which contain the present model of interest as a particular case.

The first main result of the present paper (Theorem \ref{theo:2} in Section \ref{sec:24}) shows that for our highly degenerate and time inhomogeneous $5$-dimensional stochastic system, the weak H\"ormander condition holds at all points of the state space. As a consequence, continuous transition densities exist with respect to the $5$-dimensional Lebesgue measure, at every time $t\ge 0$ and for arbitrary deterministic starting points for the system. This strong result holds provided that the deterministic signal, which is perturbed by an Ornstein-Uhlenbeck process, is an analytic function of time. With this assumption, the system we are considering has analytic coefficients. For general systems as considered in H\"opfner, L\"ocherbach and Thieullen 2013 we can achieve the weak H\"ormander condition only locally. 

Our second main result (Theorems \ref{theo:ergodicOU}+\ref{theo:ergodiccont} in Section \ref{sec:25}) deals with the long-time behavior of the process and shows that the process possesses ergodic invariant measures all of which admit a continuous density with respect to the Lebesgue measure. Moreover, there exist at most finitely many extremal invariant measures, all supported by disjoined sets. 

Our results are stated in terms of Harris recurrence which we formulate either in terms of the $T-$skeleton chain $(X_{kT})_{k\in\bbn_0}$ (the process observed at multiples of the periodicity $T$) or in terms of the $6$-dimensional continuous-time process $(i_T(t),X_t)_{t\ge 0}$ where $i_T(t)$ is $t$ modulo $T.$ We recall that a strong Markov process is called `recurrent in the sense of Harris' if it possesses an invariant measure $m$ such that any set $A $ with $m(A) > 0 $ is visited infinitely often by the process, almost surely. Here the process is allowed to start from any possible deterministic initial point. In particular, Harris recurrence implies irreducibility. If the process is recurrent in the sense of Harris, then the invariant measure is unique (up to multiplication with a constant).  Recurrence in the sense of Harris is a powerful tool in the study of the long-time behavior of a process; positive Harris recurrence implies the ergodic theorem, which is an important step towards the implementation of statistical procedures in order to identify relevant unknown parameters of the underlying model.

For the $T$-skeleton chain, we prove the existence of a finite number of disjoint Harris sets (more precisely: there is at least one Harris set, and at most a finite number) in the sense of Meyn and Tweedie 1992, Theorems 2.1 and 4.5. In restriction to any of these Harris sets, the skeleton chain is recurrent in the sense of Harris, and we have one extremal invariant measure on each Harris set. Similarly in continuous time, the $6$-dimensional system $(i_T(t),X_t)_{t\ge 0}$ admits a finite number (at least one) of disjoint invariant control sets in the sense of Arnold and Kliemann 1987; we have one extremal invariant measure on each invariant control set, and in restriction to any of the control sets, the process is recurrent in the sense of Harris. The finitely many disjoint Harris and/or invariant control sets represent a finite number of typical `stochastic equilibrium settings' for the process (in a sense of invariant law, in a sense of long time behavior), in contrast to the deterministic situation where infinitely many equilibrium states coexist.

The fact that Ornstein-Uhlenbeck diffusion has analytical coefficients comes in at several key steps of our proofs, and control arguments together with the support theorem for diffusions (see e.g.\ Millet and Sanz-Sole 1994) play a main role. 

Approaches which view neurons as deterministic (e.g.\ 
Guckenheimer and Oliva 2002, Rubin and Wechselberger 2007, Desroches, Guckenheimer, Krauskopf, Kuehn, Osinga and Wechselberger 2012) or close-to-deterministic dynamical systems (e.g.\ Berglund and Gentz 2010, Berglund and Landon 2012) have received a lot of attention, and quite often -- because of the analytical complexity of the deterministic Hodgkin-Huxley model -- one is forced to switch to simplified systems of equations such as the FitzHugh-Nagumo model or the Morris-Lecar model whose dynamics are tractable, at the price of questionable biological relevance. In this approach, `noise'
added to the classical deterministic dynamical system is often considered as `small' in order to make the stochastic system mimic essential features of the deterministic system. In contrast to this aim, in our approach `noise' -- in the form of the one-dimensional Brownian motion driving the Ornstein-Uhlenbeck SDE and by means of this the $5$-dimensional stochastic system -- is strong enough to smoothen the stochastic dynamics, despite the degeneracy of the system, by the interaction between drift and diffusion through  its $5$ dimensions. The Harris properties which we prove open a road which allows us to work in the restriction to Harris sets with ratio limit theorems or with limit theorems, and to deal in a genuinely stochastic way with long-time properties of a stochastic Hodgkin-Huxley model.  

One of our main results proves the existence of only finitely many extremal invariant measures for the stochastic Hodgkin-Huxley system. On the contrary to this situation, the deterministic Hodgkin-Huxley system exhibits a broad range of possible and qualitatively quite different behavior of its solution, depending on the specific form of the input (time-constant input, time-periodic input, jump functions ...), and depending on the starting point. Desired periodic behavior (which resembles to spiking patterns observed in neurons) appears only in special situations. Rinzel and Miller 1980 specified some interval $I$ such that a time-constant input $c\in I$ results in periodic behavior of the solution. Aihara, Matsumoto and Ikegaya 1984 determined some interval $J$ such that an oscillating input $t\to S(f\,t)$ with frequencies $f\in J$ (for some given continuous $1$-periodic function $S$) yields periodic behavior of the solution. Periodic behavior includes cases where one period of the output is equal to some multiple of the period of the input. Both papers also specify intervals $\wt I$ and $\wt J$ such that a time-constant input $c\in \wt I$ or an oscillating input at frequency $f\in \wt J$ leads to a chaotic behavior of the solution. The intricate structure of the tableau of possible types of behavior (with `modern' model constants as given in Izhikevich 2007 p.\ 37--38 to be used below, slightly different from Hodgkin and Huxley's original ones) was checked by numerical calculations in Endler 2012, Chapter 2, who obtained interesting schemes of classification. In contrast to the deterministic situation, our results show that `noise smoothens the tableau' and simplifies it in an essential way: In our stochastic Hodgkin-Huxley model, the finite number of Harris sets for the skeleton chain $(X_{kT})_{k\in\bbn_0}$ or the finite number of invariant control sets for the $6$-dimensional continuous-time process $(i_T(t),X_t)_{t\ge 0}$ corresponds to a finite number of possibilities for `typical' long time behavior. 

Our paper is organized as follows. We present the deterministic Hodgkin-Huxley system and our stochastic model in Sections \ref{sec:1.1}--\ref{sec:2.5}. Sections \ref{sec:2.3}--\ref{sec:24} focus on the weak H\"ormander condition. First, Theorem \ref{theo:hoer} formulates a sufficient condition (considering  Lie brackets of some fixed order) for validity of the weak H\"ormander condition: on the state space of the $5$-dimensional stochastic Hodgkin-Huxley system, this condition holds up to at most an exceptional set of Lebesgue measure zero. Theorem \ref{theo:2} then strengthens this and proves (via a control argument) that in fact the exceptional set is void and the weak H\"ormander condition holds everywhere. This is our first major result for the stochastic Hodgkin-Huxley system. As a consequence, Corollary \ref{cor:density} states continuity properties of Lebesgue densities of transition probabilities. 
Section \ref{sec:25} deals with ergodicity properties of the system. Thanks to a Lyapunov function we show that some large compact set is visited infinitely often. Then, using Nummelin splitting based on the results of Corollary \ref{cor:density}, we can cover the compact set with a finite number of balls of a certain type which induce renewal times in the sense of Nummelin 1978. 
Theorem \ref{theo:ergodicOU} then establishes Harris recurrence in restriction to a finite number of Harris sets for the skeleton chain $(X_{kT})_{k\in\bbn_0}.$ Theorem \ref{theo:ergodiccont} formulates this for the continuous-time process $(i_T(t),X_t)_{t\ge 0}$ in restriction to invariant control sets.  Longer proofs are shifted to the Sections \ref{section:th1}--\ref{sec:arnold-kliemann}: Section \ref{section:th1} calculates Lie brackets, Sections \ref{sec:supp} and \ref{sec:proofOU} work with control systems and the support theorem to check which parts of the state space are attainable for the stochastic Hodgkin-Huxley system. Finally, Section \ref{sec:arnold-kliemann} deals with invariant control sets for the process $(i_T(t),X_t)_{t\ge 0}$ in order to establish the link between invariant measures for the $5$-dimensional skeleton chain $(X_{kT})_{k\in\bbn_0}$ and the $6$-dimensional continuous-time process $(i_T(t),X_t)_{t\ge 0}$.

\section{Deterministic and stochastic Hodgkin-Huxley system. Main results.}\label{sec:main}

We consider a neuron modeled by a Hodgkin-Huxley system which receives a periodic input $S$ from its dendritic system. The input is random and, as a function of time, modeled by a time inhomogeneous diffusion of mean reverting type, as argued by H\"opfner 2007. We start by recalling briefly the classical deterministic Hodgkin-Huxley model.

\subsection{HH with deterministic $T$-periodic input}\label{sec:1.1} 

The classical Hodgkin-Huxley model is a $4$-dimensional ordinary differential equation. The first variable $V$ represents the membrane potential, while the other three variables $n,m$ and $h$ are related to the proportion of different types of open ion channels, which allow sodium or potassium ions to enter or to leave the neuron. 

Let $t\to S(t)$ be a $T$-periodic deterministic signal. The Hodgkin-Huxley equations with input $S(t) $ are 
$$
\left\{\begin{array}{l}
dV_t  \;=\; S(t)\, dt \;- \left[\, \ov g_{\rm K}\,n_t^4\, (V_t-E_{\rm K}) \;+\; \ov g_{\rm Na}\,m_t^3\, h_t\, (V_t  -E_{\rm Na}) \;+\; \ov g_{\rm L}\, (V_t-E_{\rm L}) \right] dt\\
dn_t \;=\;  \left[\, \al_n(V_t)\,(1-n_t)  \;-\; \beta_n(V_t)\, n_t  \,\right] dt  \\
dm_t \;=\;  \left[\, \al_m(V_t)\,(1-m_t)  \;-\; \beta_m(V_t)\, m_t  \,\right] dt  \\
dh_t \;=\;  \left[\, \al_h(V_t)\,(1-h_t)  \;-\; \beta_h(V_t)\, h_t  \,\right] dt ,
\end{array}\right. 
\leqno{\rm (HH)}$$
where
$$\ov g_{\rm K} = 36, \; \ov g_{\rm Na} = 120, \; \ov g_{\rm L} = 0. 3 , \;E_{\rm K} = - 12, \; E_{\rm Na} = 120, \; E_{\rm L} =  10.6 ,$$ 
with notations and constants of Izhikevich 2009, pp.\ 37--38. The functions $\al_n, \beta_n, \al_m, \beta_m, \al_h, \beta_h$ in (HH) take values in $(0,\infty)$ and are analytic, i.e. they admit a power series representation on $\bbr$. They are given as follows.
\begin{equation}
\begin{array}{llllll}
\alpha_n(v)  &=& \frac{0.1-0.01v }{\exp(1-0.1v)-1}, &  \beta_n(v) &= &0.125\exp(-v/80) ,  \\
\alpha_m(v)& = &\frac{2.5-0.1v}{\exp(2.5-0.1v)-1} , & \beta_m(v)&= &4\exp(-v/18) ,  \\
\alpha_h(v) &= &0.07\exp(-v/20) , &\beta_h(v) &=& \frac{1}{\exp(3-0.1v)+1}.
\end{array}
\end{equation}
For $v\in\bbr ,$ let
\begin{equation}\label{eq:ninfty}
n_\infty(v) := \frac{\al_n}{\al_n+\beta_n}(v) \;,\; m_\infty(v) := \frac{\al_m}{\al_m+\beta_m}(v) \;,\; h_\infty(v) := \frac{\al_h}{\al_h+\beta_h}(v) \;.
\end{equation}
If we think of keeping the variable $V$ constant in (HH), then these are equilibrium values in $(0,1)$ for the variables $n$, $m$, $h$ when $V\equiv v\in\bbr$.  

We write 
\beao
E_4 := \bbr\times [0,1]^3 \quad\mbox{for the state space of $(V, n,m,h)$, with points $(v,n,m,h)$} 
\eeao  
(see Proposition \ref{prop:non-exploding} below for a proof of the fact that the system stays in $E_4$ whenever it starts there).
We use the notation $F:E_4\to \bbr$ for drift terms not related to the signal in the first equation of (HH): 
\beam\label{eq:F}
F(v,n,m,h) &:=& \ov g_{\rm K}\,n^4\, (v-E_{\rm K}) \;+\; \ov g_{\rm Na}\,m^3\, h\, (v  -E_{\rm Na}) \;+\; \ov g_{\rm L}\, (v-E_{\rm L}) \nonumber \\
&=& 36\,n^4\, (v+12) \;+120\,m^3\, h\, (v -120) \;+ 0.3\, (v-10.6) \;. 
\eeam

Define from (\ref{eq:F}) a function $F_\infty : \bbr\to\bbr$ by
\begin{equation}\label{eq:Finfty}
F_\infty(v) \;:=\;  F\left( v, n_\infty(v), m_\infty(v), h_\infty(v) \right) \;,\; v\in\bbr \;.
\end{equation}
In particular, if we select $c\in\bbr$ such that  $c = F_\infty(v),$ then
\begin{equation}\label{eq:equi}
( v, n_\infty(v), m_\infty(v), h_\infty(v) )\;\in\; E_4
\end{equation}
is an equilibrium point for the deterministic system (HH) with constant signal $S(\cdot)\equiv c  .$

\begin{ex}\label{ex:2}
It is well known that for sufficiently large values of a constant signal $S(\cdot)=c$, the deterministic system (HH) exhibits regular spiking (see Rinzel and Miller 1980, for the model constants used here see Endler 2012, Section 2.1, in particular Figure 2.6). This means that for such values of $c$, the equilibrium point (\ref{eq:equi}) is unstable, and that there is a stable orbit for the $4$-dimensional system $(V, n,m,h)$ of `biological variables'.
\end{ex}

\subsection{Ornstein-Uhlenbeck type $T$-periodic diffusions and HH system with stochastic input}\label{sec:2.5}
From now on we suppose that the $T$-periodic signal $t\to S(t)$ of Subsection \ref{sec:1.1} is an analytic function. We consider a stochastic Hodgkin-Huxley system which receives this signal from its dendritic system as random input. This random input is modeled by the following diffusion 
\begin{equation}\label{eq:xi}
d\xi_t \;=\; (\, S(t)-\xi_t\,)\, \tau dt \;+\; \gamma \, \sqrt{\tau} dW_t ,  
\end{equation}
where we have chosen a parametrization in terms of $\tau > 0 $ (governing the `speed' of the diffusion) and $\gamma > 0 $ (governing the `spread' of the one-dimensional marginals). 
The process $\xi$ is a time inhomogenous Ornstein-Uhlenbeck type diffusion which carries the signal $S$.

\begin{rem}
We have an explicit representation 
$$ \xi_t = x e^{ - \tau ( t- s)} + \int_s^t e^{ - \tau (t-v) }
\left( \tau S(v) dv + \gamma \sqrt{\tau} d W_v \right) , \; t \geq s ,
$$
for the process starting at time $s$ in $x.$ Introducing the function $ s 
\to M(s) = \int_0^\infty S( s - \frac{r}{\tau} ) e^{ - r } dr ,$ the invariant law $\pi $ of the skeleton chain
$(\xi_{kT})_{k \in \bbn} $ is 
$$ \pi = {\cal N} ( M(0), \frac{\gamma^2 }{2} ) $$
and the law of $\xi_s $ starting at time $t=0$ from $\xi_0 \sim \pi$ is 
$$ {\cal L}_{ \pi , 0} ( \xi_s) = {\cal N} ( M(s) , \frac{ \gamma^2 }{2} ) $$
(cf. H\"opfner and Kutoyants 2010, Ex. 2.3). Hence the $T-$periodic signal $S(\cdot) $ is 
expressed in the process $\xi $ under the `periodically invariant' regime in the form 
of moving averages 
$$ s \to E_{\pi, 0} ( \xi_s ) = M(s) = \int_0^\infty S(s - \frac{r}{\tau} ) e^{ - r } dr $$
which are $T-$periodic. For large values of $\tau$, $M(\cdot)$ is close to $S(\cdot)$. 
\end{rem}

Consider now the HH equations driven by stochastic input $d\xi_t$, i.e.\ the $5$-dimensional system 
$$
\left\{\begin{array}{l}
dV_t  \;=\; d\xi_t \;- \left[\, \ov g_{\rm K}\,n_t^4\, (V_t-E_{\rm K}) \;+\; \ov g_{\rm Na}\,m_t^3\, h_t\, (V_t  -E_{\rm Na}) \;+\; \ov g_{\rm L}\, (V_t-E_{\rm L}) \right] dt\\
dn_t \;=\;  \left[\, \al_n(V_t)\,(1-n_t)  \;-\; \beta_n(V_t)\, n_t  \,\right] dt  \\
dm_t \;=\;  \left[\, \al_m(V_t)\,(1-m_t)  \;-\; \beta_m(V_t)\, m_t  \,\right] dt  \\
dh_t \;=\;  \left[\, \al_h(V_t)\,(1-h_t)  \;-\; \beta_h(V_t)\, h_t  \,\right] dt  \\
d\xi_t \;=\; (\, S(t)-\xi_t\,)\, \tau dt \;+\; \gamma\, \, \sqrt{\tau} dW_t .
\end{array}\right. 
\leqno{\rm (\xi HH)} 
$$
Write $X = (X_t)_{t \geq 0} , $ $ X_t = ( V_t, n_t, m_t , h_t, \xi_t ) ,$ for the solution of ${\rm (\xi HH)}$ (we show in Proposition \ref{prop:non-exploding} below the existence of a unique strong solution),  $E_5=\bbr\times [0,1]^3\times \bbr$ for the corresponding state space, and denote the elements of $E_5$ by $x= (v,n,m,h,\zeta)$. We write $P_x$ for the probability measure under which the solution $X = (X_t)_{t \geq 0 } $ of ${\rm (\xi HH)}$ starts from $x.$ Let $\left(\,P_{s_1,s_2}(x_1, dx_2)\,\right)_{0\le s_1<s_2<\infty}$ be the associated semigroup of transition probabilities. Due to the $T$-periodicity of the deterministic signal $t\to S(t)$, the semigroup is $T$-periodic in the following sense:  
$$
P_{s_1,s_2}(x_1, dx_2) \;=\; P_{s_1+kT,s_2+kT}(x_1, dx_2) \quad\mbox{for all $k\in\bbn_0$} \;. 
$$

\begin{prop}\label{prop:non-exploding}
For any $x \in E_5,$ there exists a unique strong non-exploding solution $X$ to ${\rm (\xi HH)} $ starting from $x$ at time $0 ,$ taking values in $E_5.$  
\end{prop}

{\bf Proof.} By our assumptions, a strong solution $\xi_t$ of (\ref{eq:xi}) exists. Moreover, the coefficients of $V$ and $n,m, h$ are locally Lipschitz continuous. This implies the existence of a unique strong solution of the system ($\xi \! \!$ HH) which is a maximal solution, i.e. it exists up to some explosion time. So we have to prove that the process does not explode. By assumption, $\xi_t$ does not explode. Consider now the unique solution $(V_t, n_t, m_t, h_t, \xi_t) $ of ${\rm (\xi HH)} $ on $ [ 0, T_\infty [, $ where $T_\infty $ is the associated explosion time. It is easy to show 
that $n,m$ and $h$ stay in $ [0,1] ,$ whenever they start in $[0, 1] ,$ see for instance Proposition 1 of H\"opfner, L\"ocherbach and Thieullen 2013. 

Replacing $n,m$ and $h$ by their maximal value $1,$ we obtain easily 
from ${\rm (\xi HH)} $ that
\begin{equation}\label{eq:vtbounded}
 | V_t | \le | \xi_t | + C_1 \int_0^t | V_s| ds + C_2 t , \; t  < T_\infty ,
\end{equation}
where $C_1 $ and $C_2$ are suitable constants.
This implies, using Gronwall's inequality and the non-explosion of $\xi_t,$ that $V_t$ does not explode either. 
Hence $T_\infty = \infty $ almost surely and the above estimates hold on $ [ 0, \infty [ .$ 
\halmos

\subsection{Weak H\"ormander condition}\label{sec:2.3} 
Our system ${\rm (\xi HH)} $ is a $5$-dimensional diffusion driven by one-dimensional Brownian motion. As a consequence, the only possibility for guaranteeing non-degeneracy of the system is that the system `feels the noise via the drift'. In other words, we have to check whether the weak H\"ormander condition holds. Since the drift term of ${\rm (\xi HH)} $ depends on time, we add time as a first coordinate to our system. More precisely, we write $\bbt := [0,T]$ for the torus and identify $t$ with $i_T(t) := t \mod T.$ Elements of $ [0, \bbt  ] \times E_5 $ will be denoted either by $ (x^0, x^1 , \ldots , x^5 ) $ or by $ ( t, x)$ or by $ ( t, v, n, m,h, \zeta ) .$ Working with the space-time process $\bar X_t = (i_T (t), X_t), $ the associated drift and diffusion coefficients are the vector fields
\begin{equation}\label{eq:b} \bar b (t,x)= \left( 
\begin{array}{c}
1\\
b^1 (t,x) \\
\vdots \\
b^5 (t,x) 
\end{array}
\right) \in \bbr^6  \quad \mbox{ and } \quad \bar \sigma (t,x)= \gamma \sqrt{\tau}  \left( 
\begin{array}{c}
0\\
1\\
0\\
0\\
0\\
1
\end{array}
\right) \in \bbr^6 ,
\end{equation}
where for $ x = ( v, n, m,h, \zeta),$
$$ b^1 (t,x) = (S(t) - \zeta ) \tau  - F( v , n , m, h) ,\; \,  b^2 (t, x) =\al_n(v)\,(1-n)  - \beta_n(v)\, n, $$
$$  b^3 ( t, x) = \al_m(v)\,(1-m)  - \beta_m(v)\, m, \;  b^4 (t, x) = \al_h(v)\,(1-h)  - \beta_h(v)\, h , \;  b^5 (t,x) = (S(t) - \zeta ) \tau .$$
 
We identify $ \bar b(t, x) $ and $ \bar \sigma ( t, x) $ with differential operators
$$ \bar b(t,x) =  \frac{\partial }{\partial t } + \sum_{ i = 1}^5  b^i (t,x) \frac{\partial }{\partial x^i } \quad
\mbox{and} \quad  \bar \sigma ( t, x) = \sum_{ i = 1}^5  \sigma^i (t,x) \frac{\partial }{\partial x^i } . $$  

We are now going to introduce the successive Lie brackets that we use in the sequel. First, recall that for 
vector fields $f(t, x)$ and $ g(t,x) : \bbt \times E_5   \to \bbr^6 , $ the Lie bracket $ [f, g ] $ is defined by
$$ [f, g ]^i = \sum_{j = 0}^5 \left( f^j \frac{ \partial g^i }{\partial x^j } - g^j \frac{ \partial f^i }{\partial x^j } \right) \; , \; i = 0, \ldots , 5, $$
with superscript `$i$' for the $i$-th component.
In this way, for a vector field $f: \bbt \times E_5 \to \bbr^6$ whose `0-component' equals $0$, the Lie bracket $[\bar b, f]$ takes the form
$$
[\bar b, f]^0 = 0 \quad,\quad 
[\bar b,f]^i = \frac{\partial f^i}{\partial t} + \sum_{j=1}^5 \left( b^j\, \frac{\partial f^i}{\partial x^j} - f^j\, \frac{\partial b^i}{\partial x^j}\right) \;,\; i=1,...,5 ,
$$
and the Lie bracket  $[\bar \sigma ,f]$ takes the form  
$$
[\bar \sigma ,f]^0 = 0 \quad,\quad 
[\bar \sigma ,f]^i =  \gamma \sqrt{\tau} \left( \frac{\partial f^i }{\partial x^1 } + \frac{\partial f^i }{\partial x^5 } \right)    \;,\; i=1,...,5 \;. 
$$

We introduce the following system of sets of vector fields based on iterated Lie brackets.

\begin{defin}\label{def:hoer}
Define a set $\call$ of vector fields by the `initial condition' $\bar \sigma \in \call$ and  an arbitrary number of iterations steps 
\beqq\label{eq:iteration}
L\in\call \;\Lra\; [\bar b,L] , [\bar \sigma , L] \in \call  \;.   
\eeqq
For $N\in\bbn$, define the subset $\call_N$ by the same initial condition and at most $N$ iterations (\ref{eq:iteration}). Write $\call_N^*$ for the closure of $\call_N$ under Lie brackets; finally, write 
$$
\Delta_{\call_N^*}\;:=\;\mbox{\rm LA}(\call_N)
$$
for the linear hull of $\call_N^*$, i.e.\ the Lie algebra spanned by $\call_N$. 
\end{defin}

Note that all elements of $\call_N^*$ have `$0$-component' equal to zero, so $5$ is an obvious upper bound for ${\rm dim}(\Delta_{\call_N^*}) .$ 

\begin{defin}
We say that a point $x^* \in E_5$ is of {\em full weak H\"ormander dimension} if there is some $N\in\bbn$ such that   
\beqq\label{fwHd}
({\rm dim}\; \Delta_{\call^*_N})(s,x^*) \; =\; 5 \quad\mbox{independently of $s\in\bbt$} \;. 
\eeqq
We put 
$$ {\cal I}_5 := \{ x = (v, n, m, h , \zeta ) \in E_5 : \mbox{ $x$ is of full weak  H\"ormander dimension }  \} .$$ 
\end{defin}

\begin{rem}
Notice that in the iteration step (\ref{eq:iteration}), it is allowed to build Lie brackets using the drift vector $\bar b (t, x) .$ It is for this reason that the above condition is called `weak' in contrast to the `strong' H\"ormander condition. In the strong H\"ormander condition, only iterations using the column vectors of the diffusion matrix are allowed. Since in our case the diffusion matrix is built of only one column, it is clear that the strong H\"ormander condition can never hold.
\end{rem}

In the following we are going to state a sufficient condition ensuring that a given point belongs to ${\cal I}_5.$ In order to do so, let
\begin{equation}\label{eq:det}
D(v,n,m,h) \;:=\; \det \left( \begin{array}{lll} 
\partial^2_v b^2   & \partial^3_{v} b^2 & \partial^4_{v}  b^2 \\
\partial^2_v b^3  &  \partial^3_{v} b^3 &  \partial^4_{v} b^3  \\
\partial^2_v b^4 &  \partial^3_{v}  b^4 & \partial^4_{v}  b^4  \\
\end{array} \right) (v,n,m,h) \quad,\quad (v,n,m,h)\in E_4 \;, 
\end{equation}
where $\partial^k_v$ denotes the $k-$fold partial derivative with respect to $v $\footnote{Notice that $ D( v, n, m,h ) $ does not depend on time.}, and introduce
$$ {\cal O} := \{ (v,n,m,h) \in E_4 : D( v,n,m,h) \neq 0 \} . 
$$

We quote the following proposition from H\"opfner, L\"ocherbach and Thieullen 2013.
\begin{prop}[Proposition 9 of \cite{h-l-t}]
${\cal O}$ is an open set of full Lebesgue measure, i.e. $\lambda (  {\cal O}^c) = 0 .$ 
\end{prop}


Calculating the first four Lie brackets of our system by using successively first the drift vector and then three times the diffusion coefficient, we obtain the following theorem. 
 
\begin{theo}\label{theo:hoer}
All points $ x = ( v,n,m,h, \zeta) $ in $E_5$ whose  first four components belong to ${ \cal O} $ are points satisfying the weak H\"ormander condition. 
\end{theo}

The proof of Theorem \ref{theo:hoer} is given in Section \ref{section:th1}.

\begin{rem}\label{rem:numeric}
We resume the numerical study of Section 5.4 in H\"opfner, L\"ocherbach and Thieullen 2013. First of all, the set $\calo$ is certainly non-empty since we find a strictly negative value of the determinant e.g.\ at the equilibrium point $(0, n_\infty(0), m_\infty(0), h_\infty(0) )$ of the 4d deterministic system (HH). In order to obtain more information about $\calo$, we calculate 
$$
v \;\lra\; D\left( v, n_\infty(v), m_\infty(v), h_\infty(v) \right)   
\leqno(*)
$$ 
at equilibrium points of (HH) which correspond to constant input $S(\cdot)\equiv c$. By equations (\ref{eq:Finfty}) and (\ref{eq:equi}), equilibrium requires $c:=F_\infty(v)$. Since $v \to F_\infty(v)$ is strictly increasing, there is a one-to-one correspondence between values $v$ of the membrane potential and values $c$ of the input: as an example, to $v\in\{-10,0,10\}$ correspond $c\in\{-6.15, -0.05, 26.61\}$.

Calculating the function $(*)$ when $v$ ranges over the interval $(-15,30)$, we find zeros at the two points $v=-11.48$ and $v=10.34$ (numerical values rounded to 2 decimals). The function  $(*)$ is strictly negative in between, and strictly positive outside. In particular, equilibrium points of (HH) under  constant input $-6.15\le c\le 26.61$ belong to the set $\calo$.
\end{rem}

\subsection{Non-degeneracy of the stochastic Hodgkin-Huxley process}\label{sec:24}
Let $ y_* = (0, n_\infty(0), m_\infty(0), h_\infty(0) )  $ be the equilibrium point for the deterministic system (HH) driven by constant input $ c = F_\infty (0) \approx - 0. 0534 .$ By Remark \ref{rem:numeric} above we know that $ y_*  \in {\cal O}.$ In this section we will show that for any neighborhood $U$ of $y_* ,$ the set $ U \times \bbr $ is accessible. Since the coefficients of the system are analytic, this will imply that the weak H\"ormander condition holds on the whole state space $E_5.$ 

We start with the following proposition which is due to discussions  with Michel Bena\"im, see also Bena\"im, Le Borgne, Malrieu, Zitt 2012. It shows that, starting from any initial point, our system can reach $U \times \bbr $ for any open neighborhood $U$ of $y_* .$
\begin{prop}\label{theo:accessible}
Let $U \subset E_4 $ be a neighborhood of  $  y_* .$ Then for all $x \in E_5,$ there exists $t_0 $ such that for all $ t \geq t_0  $ 
$$ P_{0, t } ( x , U \times \bbr ) > 0 .$$

In particular, for the $T-$skeleton chain $ (X_{ k T })_{ k \geq 0} $ it holds that for all $x \in E_5 $ there exists $k \geq 1 $ such that 
$$ P_x ( X_{kT } \in U \times \bbr ) > 0 .$$
\end{prop} 

The proof of this proposition is given in Section \ref{sec:supp} below. As a consequence, we obtain the following theorem.  

\begin{theo}\label{theo:2}
The weak H\"ormander condition holds on $E_5 ,$ i.e. ${\cal I}_5 = E_5 .$ 
\end{theo} 

{\bf Proof.}
The proof uses the following fact. For any diffusion process having analytic coefficients, the following holds true:
\begin{equation}\label{eq:firststatement}
\mbox{if $X_t \notin {\cal I}_5 ,$ then $ X_{t +s } \notin {\cal I}_5 $ for all $s \geq 0 $ almost surely . }
\end{equation}
For the convenience of the reader we will give a proof of (\ref{eq:firststatement}) 
in Section \ref{sec:proofOU} below.

Based on (\ref{eq:firststatement}), we argue as follows. Suppose that there exists $x \in E_5 \setminus {\cal I}_5 .$ We will apply Proposition \ref{theo:accessible} with this fixed starting point. Since $y_* \in {\cal O}, $ we may choose a neighborhood $ U$ of $ y_* $ sufficiently small such that $ U \subset {\cal O} .$ Since $ U \times \bbr \subset {\cal O} \times \bbr \subset {\cal I}_5,$ Proposition \ref{theo:accessible} then implies that there exists $t^* $ such that for the fixed $ x \in E_5 \setminus {\cal I}_5,$
$$ P_x ( X_{t^*} \in {\cal I}_5  ) \geq P_x ( X_{t^*} \in U \times \bbr ) > 0 .$$ 
But, applying (\ref{eq:firststatement}), we have $ P_x ( X_{t^*} \in {\cal I}_5 ) = 0 ,$ since $ x \notin {\cal I}_5.$ This is a contradiction. \halmos 

Once the weak H\"ormander condition holds on $E_5,$ it follows that the process possesses Lebesgue densities. 

\begin{cor}\label{cor:density}
For $0\le s_1<s_2<\infty$, consider the process $X$ starting at time $s_1\ge 0$ from arbitrary $x\in E_5$. Then the law $\,P_{s_1,s_2}(x, \cdot)\,$ admits a Lebesgue density $\,p_{s_1,s_2}(x, y)\,$. For fixed $x,$ $p_{s_1,s_2}(x, y)$ is continuous in $y,$ uniformly in $x.$ Moreover, for any fixed $ x' \in E_5  , $ the map $ x \to p_{s_1,s_2}(x, x') $ is lower semi-continuous.
\end{cor}

{\bf Proof.} The weak H\"ormander condition holds everywhere. If the coefficients of our system where $C^\infty_b $ and time homogeneous, then classical results as presented e.g. in Kusuoka and Stroock 1985, Corollary (3.25), or in Nualart 1995, Theorem 2.3.3, would allow us to conclude. However, the coefficients of our system are not $C^\infty_b $ and they are time inhomogeneous. But treating time as a first coordinate and using a localization argument allows to prove the assertion. The proof of Theorem 1 in H\"opfner, L\"ocherbach and Thieullen 2013 gives the details. 
\halmos

\subsection{Ergodicity of the stochastic Hodgkin-Huxley system}\label{sec:25}
We start by showing, using Lyapunov functions, that almost surely the system comes back to a compact set infinitely often. We are working with the $T-$skeleton $(X_{kT})_{ k \geq 0 },$ where $T$ is the periodicity of the underlying signal $S.$ 

\begin{prop}
\begin{enumerate}
\item
There exists a compact set $K \subset E_5 $ such that for all $x \in E_5, $ $P_x-$almost surely, 
$$ \sum_{k= 0}^\infty 1_K ( X_{k T }) = \infty .$$ 
\item
There exist an integer $N \geq 1 , $ $\varepsilon_1 > 0 , \ldots , \varepsilon_N > 0 ,$  $ x_1, \ldots , x_N \in K $ and $ y_1 , \ldots , y_N \in E_5$ such that 
$K $ is covered by $ B_{\varepsilon_1} ( x_1 ) , \ldots , B_{\varepsilon_N } (x_N) $ and such that for any $ 1 \le i \le N,$
$$  \inf_{ x' \in B_{\varepsilon_i} (x_i) , y' \in B_{\varepsilon_i }( y_i ) } p_{ 0, T } (x', y') > 0 .$$
\end{enumerate}
\end{prop}

{\bf Proof.}  Let $\Phi :  E_5  \to [1, \infty [ $ be a $C^2 -$function satisfying $\Phi ( x) = |x^1 | + (x^5)^2 $ for all $x$ such that $ |x^1| \geq 2 ,$  $\Phi (x) $ arbitrary elsewhere. Write $L_t$ for the generator of ($\xi$HH) at fixed time $t.$ Since $ 0\le x^2, x^3, x^4 \le 1 , $ it is easy to see that  
\begin{equation}\label{eq:lyap}
 L_t \Phi ( x) \le - c_1 \Phi (x) + c_2 1_{\tilde K} ( x) ,
\end{equation}
where $c_1, c_2 $ are positive constants and where $\tilde K = \{ x \in  E_5  : |x^1| \le C , |x^5 | \le C  \} $ is a compact subset of $  E_5  .$ Applying It\^o's formula to $ e^{ c_1 t } \Phi ( x) $ and using localization with $ \inf \{ t : \Phi ( X_t ) \geq m \} $ as $m \to \infty, $ we obtain
$$ E_x  \Phi ( X_t) \le e^{  - c_1 t}   \Phi ( x) + \frac{c_2}{c_1} $$
for all $ t > 0.$ Let now $t = T, $ where $T$ is the period of the underlying signal. Thus
$$ P_{ 0, T } \Phi ( x) - \Phi ( x) \le - ( 1 - e^{ - c_1 T } ) \Phi ( x)+ \frac{c_2}{c_1} .$$  
In particular, there exists a constant $ C_2 ,$ such that 
\begin{equation}\label{eq:lyapunovgrid}
 P_{ 0, T } \Phi ( x) - \Phi ( x) \le - \varepsilon \quad \mbox{ for all $ x $ with $ |x^1| > C_2$ or $ |x^5| > C_2 ,$} 
\end{equation}
for some fixed $ \varepsilon > 0 .$ By Theorem 4.3 of Meyn and Tweedie 1992, we know that (\ref{eq:lyapunovgrid}) implies the following statement: Starting from any point in $E_5,$  the skeleton chain $ (X_{kT})_{ k \in \bbn}$ visits the compact set $K = \{ x \in  E_5  : |x^1| \le C_2, |x^5| \le C_2 \}  $ infinitely often. This proves the first assertion of the proposition. 

Concerning the second assertion, observe that for any $x \in K  $ there exists $y \in E_5$ such that $p_{0, T } ( x, y) > 0.$ By continuity in $y$ and lower semi-continuity in $x,$ this can be extended to small balls around $x$ and $y.$

As a consequence, for any $x \in K$ there exist $ y $ and $\varepsilon > 0 $ such that 
\begin{equation}\label{eq:min}
 \inf_{ x' \in B_{ \varepsilon} ( x), y' \in  B_{ \varepsilon} ( y)} p_{ 0, T } ( x', y') > 0 . 
\end{equation}

Hence the compact set $K$ is covered by a finite number of such balls $ B_{\varepsilon_1} ( x_1)  , \ldots ,  B_{\varepsilon_N} ( x_N) ,$  with associated points $y_1 , \ldots , y_N .$ This shows the second assertion of the proposition. \halmos

The lower bound (\ref{eq:min}) can be rewritten as follows. For any $ 1 \le k \le N ,$  
\begin{equation}\label{eq:lbnummelin}
 P_{ 0, T } (x , dy ) \geq \beta_k \; 1_{ B_{\varepsilon_k }(x_k)} (x)  \nu_k (dy) , 
\end{equation}
where 
$$ \beta_k = \lambda ({ B_{\varepsilon_k} ( y_k) } ) \cdot \inf_{ x' \in B_{ \varepsilon_k} ( x_k), y' \in  B_{ \varepsilon_k} ( y_k)} p_{ 0, T } ( x', y') \mbox{ and } \nu_k = \frac{1}{ \lambda ({ B_{\varepsilon_k} ( y_k) } )} \lambda_{| B_{\varepsilon_k} ( y_k) } .$$
Using Nummelin splitting (see e.g. Nummelin 1978), this implies the following.

\begin{theo}\label{theo:ergodicOU}
a) The $T-$skeleton chain $(X_{kT})_{ k \geq 0 }$ possesses ergodic invariant measures. 

b) Any invariant measure for the $T-$skeleton chain admits a continuous density with respect to Lebesgue measure on $E_5$.

c) The $T-$skeleton chain admits at most a finite number of extremal invariant measures living on disjoint Harris subsets of $E_5$. 
\end{theo}

{\bf Proof.} 
Parts a) and c) are essentially Meyn and Tweedie 1992, decomposition Theorem 2.1 and Theorem 4.5; note that our skeleton chain satisfies the assumptions of both theorems, due to assertion 2 in our Corollary \ref{cor:density}.  

a) Let $\nu$ denote any probability measure on $E_5$. Starting from $\nu$,  $P_\nu$-a.s., the skeleton chain visits the compact set $K$ infinitely often. As a consequence,  for $P_\nu$-almost all $\om$, there exists (at least one) index $k = k(\om) \in \{ 1 , \ldots , N \} $ such that the $\om$-path $ (X_{nT}(\om))_{ n \geq 0}$  of the skeleton chain  visits ${B_{\varepsilon_k } (x_k)} $ infinitely often. Let $A_k$ denote the set of all paths which visit ${B_{\varepsilon_k } (x_k)} $ infinitely often,  $1\le k\le N$. 

Let $(U_n)_n$ be an i.i.d. sequence of uniform $U(0,1)-$random variables, independent of the process. By means of these, we can  introduce a sequence of regeneration times $(1+R^{(k)}_n)_{n\ge 1}$ associated to successive visits of $B_{\vep_k}(x_k)$ in the following way: 
$$ 
R^{(k)}_{n+1} = \inf \{ l > R^{(k)}_n :  X_{ l T } \in B_{\vep_k}(x_k) , U_l \le \beta_k \} \;,\; n\ge 0 \;,\; R^{(k)}_0 \equiv 0   
$$
where $\beta_k$ is from the lower bound (\ref{eq:lbnummelin}). For every $k\in \{ 1 , \ldots , N \}$, the $R^{(k)}_n$ are finite on $A_k$ for all $n$, and satisfy $R^{(k)}_n \uparrow \infty$ on $A_k$ as $\nto$ . 

There is at least one $k \in \{ 1 , \ldots , N \} $ such that $A_k$ has positive $P_\nu$-measure. We suppose without loss of generality that $k=1$. By the lower bound (\ref{eq:lbnummelin}), $\,P_{0,T}(x,dy) \ge \beta_1\, 1_{B_{\vep_1}(x_1)}(x)\, \nu_1(dy) ,$ and using the Borel-Cantelli lemma (see also Lemma 1.1 of Meyn and Tweedie 1992), any path belonging to $A_1$ also visits $B_{\vep_1}(y_1)$ infinitely often. Recall that $\nu_1$ is the uniform measure on $B_{\vep_1}(y_1)$. By Nummelin splitting with minorization according to (\ref{eq:lbnummelin}) and regeneration times $(1+R^{(1)}_n)_n$, $\,P_\nu(A_1)>0$ thus implies
$$
P_{\nu_1}\left( R^{(1)}_1 < \infty \right) \;=\; 1  
$$
as a consequence of the Markov property at times $(1+R^{(1)}_n)_n$ and the Borel-Cantelli lemma. 

Similarly, for $1\le k\le N$, call $B_{\vep_k}(y_k)$ a `good' set if $\,P_{\nu_k}( R^{(k)}_1 < \infty ) = 1\,$.  At least one such `good' set exists, namely $B_{\vep_1}(y_1)$ in the notation above. Notice that being a `good' set is a property which only depends on the whole ball $B_{\vep_k}(y_k)$ and the semigroup of the process. Rearranging the numbering, we find some maximal subset $\{1,\ldots,N_1\}$ of $\{1,\ldots, N\}$, $1\le N_1\le N$, with the property that  $\,P_{\nu_k}( R^{(k')}_1 < \infty) = 0  \,$ whenever $k \neq k'$ in $\{1,\ldots,N_1\}$.   Stated equivalently, this rearrangement induces a partition $A_1\dot\cup A_2 \dot\cup \ldots \dot\cup A_{N_1}$ of the path space up to some remaining set of paths which has $P_\nu$-measure zero for every initial law $\nu$. Next, we define 
$$
\tau \;:=\; R^{(1)}_1 \wedge R^{(2)}_1 \wedge \ldots \wedge R^{(N_1)}_1 
$$
which is $P_\nu$-almost surely finite for every initial law $\nu$. By Nummelin splitting and the strong law of large numbers, 
\begin{equation}\label{eq:muball} 
 \mu ( f ) := \sum_{k=1}^{N_1} E_\nu \left( 1_{\{ X_{\tau T} \in B_{\vep_k}(x_k)\}}E_{X_{\tau T}} \left( 
\sum_{ l = R^{(k)}_n +1 }^{R^{(k)}_{n+1} } f ( X_{lT}) \right) \right)
\end{equation}
is an invariant measure of the skeleton chain, combining in an `adaptive' way the relevant $A_k$'s from the above partition. 
Note that this formula extends the usual form of the invariant measure to the case where several balls are present in the lower bound (\ref{eq:lbnummelin}).

Now define for any $ 1 \le k \le N_1$ the measure
$$
\mu_k (f) := E_{\nu_k}\left( \sum_{ l=0 }^{ R^{(k)}_1 } f( X_{ l T } ) \right)  
$$
and let $H_k \subset E_5$ be the support of $\mu_k, $ $1 \le k \le N_1.$
Then for any initial measure $\nu$ concentrated on $H_k$, sets $A$ with $\mu_k(A)>0$ are visited infinitely often by the skeleton chain $P_\nu$-almost surely.

b) Any invariant measure $\mu $ is absolutely continuous with respect to Lebesgue measure, thanks to Corollary \ref{cor:density},
 with Lebesgue density $ \int \mu (d z) p_{ 0, T } ( z , y ) .$  Moreover, since the continuity of $p_{0, T } ( z, y) $ in $y$ is uniform in $z,$ the corresponding Lebesgue density is continuous by dominated convergence. 

c) In order to achieve the proof of Theorem \ref{theo:ergodicOU}, 
we have to show that the skeleton chain possesses only a finite number of extremal invariant probability measures which have supports given by disjoint subsets of $E_5.$ We prove this in Section \ref{sec:arnold-kliemann} 
below. From the structure of the above Lyapunov condition (\ref{eq:lyapunovgrid}) and boundedness of $P_{0,T}\Phi$ on $K$ according to (\ref{eq:vtbounded}), we deduce that in the restriction to a Harris set, recurrence is necessarily positive recurrence.\halmos  \\

Recall that $i_T (t)$ denotes $ t \mod T $ and that $\bbt = [0,T]$ is the torus.  We get the following corollary of the above theorem.
\begin{theo}\label{theo:ergodiccont}
Under our assumptions, the process $(i_T (t),X_t)_{t\ge 0}$ 
admits at most a finite number of extremal invariant measures living on disjoint Harris subsets of $\bbt \times E_5.$  
\end{theo}

{\bf Proof: } See again Section \ref{sec:arnold-kliemann}. \\

\begin{rem}
By  Arnold and Kliemann 1987 we know that the support of any extremal invariant measure of the space-time process $(i_T (t),X_t)_{t\ge 0}$ is given by an invariant control set of the associated deterministic control system (see Section \ref{sec:arnold-kliemann} below for a precise definition, see also Colonius and Kliemann 1993). Hence the number of invariant control sets of the associated deterministic control system gives an a priori upper bound on the number of extremal invariant measures (and thus the number of Harris sets) of $(i_T (t),X_t)_{t\ge 0}.$
\end{rem}

\section{Proof of Theorem 1.}\label{section:th1}
Let $\bar X_t = (i_T (t) , V_t, n_t, m_t, h_t, \xi_t )$ be the diffusion process of ($\xi \! \!$ HH) to which we have added time as first coordinate, with state space $\bbt \times E_5.$ Recall the exact form of $\bar b(t,x)$ and $\bar \sigma (t,x)$ given in (\ref{eq:b}).
By the structure of the diffusion coefficient, the equation is already written in the Stratonovich sense.
We start by calculating the Lie-bracket of $ \bar \sigma $ and $ \bar b $ where we recall that we are working on $\bbt \times \bbr^5 ,$ including time. We have
 $$
[\bar b, \bar \sigma  ]  =  \gamma \sqrt{\tau} \left( 
\begin{array}{c}
0\\
\partial_v F \\
- \partial_v b^2\\
- \partial_v b^3\\
- \partial_v b^4\\
0
\end{array}
\right) 
+ \gamma  \tau^{3/2} \left( 
\begin{array}{c}
0\\
1 \\
0\\
0\\
0\\
1 
\end{array}
\right)
.
$$

In the same way, we obtain
$$
[ \bar \sigma ,  [ \bar b , \bar \sigma ]] = \gamma^2 \tau  \left( 
\begin{array}{c}
0\\
\partial^2_v F  \\
- \partial^2_v b^2\\
- \partial^2_v b^3\\
- \partial^2_v b^4\\ 
0 
\end{array}
\right)  \quad 
\mbox{ and } \quad  
[ \bar \sigma , [\bar \sigma , [\bar b, \bar \sigma ]]]  = \gamma^3 \tau^{3/2}  \left( 
\begin{array}{c}
0\\
\partial^3_v F(v,n,m,h) \\
- \partial^3_v b^2\\
- \partial^3_v b^3\\
- \partial^3_v b^4\\
0 
\end{array}
\right) .
$$
We obtain an analogous formula for $ [\bar \sigma  , [\bar \sigma , [\bar \sigma  , [\bar b, \bar \sigma ]]]] ,$ where fourth derivatives with respect to $v$ appear.

Now we are able to conclude our proof. By definition of $F$ in (\ref{eq:F}), $ \partial_v F (v,n,m,h) \neq 0 $ for all $ (v,n,m,h) \in E_4$ and $ \partial^k_v F(v,n,m,h) \equiv 0 $ for all $ k \geq 2 .$ Notice that the above vectors all have the first coordinate corresponding to time which equals zero. Hence we may identify them with elements of $\bbr^5 .$ Doing so, without changing notations, 
we have for all fixed $x \in E_5,$ 
$$ \det \left( 
\begin{array}{ccccc}
| & | & | & | & | \\ 
\bar \sigma &[\bar b, \bar \sigma ] &  [\bar \sigma  , [\bar b, \bar \sigma  ]] &  [\bar \sigma , [\bar \sigma  , [\bar b, \bar \sigma ]]]  & [\bar \sigma  , [\bar \sigma , [\bar \sigma  , [\bar b, \bar \sigma ]]]] \\
| & | & | & | & | 
\end{array}
\right) \neq 0 
$$ 
if and only if 
$$ \det \left( 
\begin{array}{ccccc}
1 &\partial_v F& 0 &  0 & 0 \\
0& - \partial_v b^2 & - \partial^2_v b^2 & - \partial^3_v b^2&- \partial^4_v b^2 \\
0& - \partial_v b^3 & - \partial^2_v b^3 & - \partial^3_v b^3 &- \partial^4_v b^3\\
0& - \partial_v b^4 & - \partial^2_v b^4 & - \partial^3_v b^4 & - \partial^4_v b^4 \\
1 & 0 & 0 & 0 & 0 
\end{array}
\right) \neq 0 .
$$ 
Developing this determinant first with respect to the last line and then with respect to the first line of the remaining sub-determinant, this last determinant is different from zero if and only if
$D ( v,n,m,h) \neq 0  $ (recall the definition of $D( v, n, m,h) $ in (\ref{eq:det})). As a consequence, $ \bar \sigma , [\bar b, \bar \sigma ], \ldots ,  [\bar \sigma  , [\bar \sigma , [\bar \sigma  , [\bar b, \bar \sigma ]]]] $ span $\bbr^5 $ for all $x \in E_5$ such that $D (v, n,m,h) \neq 0 ,$ and therefore, the weak H\"ormander condition is satisfied on ${\cal O}\times \bbr.$ \halmos

\section{Proof of Proposition \ref{theo:accessible}.}\label{sec:supp}
We consider the system ($\xi$HH) driven by $S$ of Section \ref{sec:2.5}, 
\begin{equation}
X_s = x + \int_0^s \sigma ( X_u)  dW_u + \int_0^s  b( u, X_u) du , \quad s\le t ,
\end{equation}
where 
\begin{equation}\label{eq:breal}  b (t,x)= \left( 
\begin{array}{c}
b^1 (t,x) \\
\vdots \\
b^5 (t,x) 
\end{array}
\right)   \quad \mbox{ and } \quad  \sigma (x)= \gamma \sqrt{\tau}  \left( 
\begin{array}{c}
1\\
0\\
0\\
0\\
1
\end{array}
\right) \in \bbr^5 .
\end{equation}
We write 
$ {\cal C}= C ( [ 0, \infty [ , \bbr^5 )$ for the space of continuous functions and endow ${\cal C} $ with its canonical filtration $ ( {\cal F}_t)_{t \geq 0 } .$ Let $ \mathbb{P}_{0,x}$ be the law of $(X_{  u }, u \geq 0 ) $ on ${\cal C} ,$ starting from $x$ at time $0.$ With $y^*$ and $U$ as in Proposition \ref{theo:accessible}, we wish to find lower bounds for quantities of the form $\mathbb{P}_{0,  x } (  B ) $ where $B = \{ f \in  {\cal C} : f(t) \in U \times \bbr \}  \in {\cal F}_t . $ In order to do so, we will use control arguments and the support theorem for diffusions. We need first to localize the system.  
Let $K_n = [ - n , n ] \times [0 ,1 ]^3 \times [ - n , n ]   \subset E_5$ and let $T_n = \inf \{ t : X_t \in K_n^c \} $ be the exit time of $K_n.$ For a fixed $n,$ let $ b^n (t, x) $ and $\sigma^n (x)$ be $C_b^\infty -$extensions in $x$ of $ b(t, \cdot_{| K_n }) $ and $\sigma_{| K_n} .$  
Let $X^n$ be the associated diffusion process. For any fixed $n_0 < n $ and any starting point $x \in K_{n_0} ,$ we write $\mathbb{P}_{0, x}^n $ for the law of $(X^n_{  u } , u \geq0  )$ on ${\cal C}. $ Then for any $t > 0 $ and for any measurable $B \in {\cal F}_t ,$ 
\begin{equation}\label{eq:tobelb}
\mathbb{P}_{0,  x } (  B ) \geq \mathbb{P}_{ 0, x } ( \{ f  \in B ; T_n > t \}  ) 
 = \mathbb{P}^n_{ 0, x } (  \{ f \in B ; T_n >  t \} )  .
\end{equation}
It suffices to show that this last expression is strictly positive, for the given set $B,$ for any fixed $x \in K_{n_0}.$ 
For this sake we will use the support theorem for diffusions of Stroock and Varadhan 1972.  Let $ {\cal H} = \{ {\tt h} : [ 0, t ] \to \bbr : {\tt h}(s) = \int_0^s \dot {\tt h} (u) du , \forall s \le t , \int_0^t \dot {\tt h}^2 (u) du < \infty \} $ be the Cameron-Martin space. 
Given ${\tt h} \in {\cal H} ,$ consider $ X({\tt h})$ the solution of the differential equation 
\begin{equation}\label{eq:control1}
X({\tt h})_s = x + \int_0^s \sigma^n ( X({\tt h})_u) \dot {\tt h} (u) du + \int_0^s  b^n ( u, X({\tt h})_u) du , \quad s\le t ,
\end{equation}
where $X({\tt h})$ is of the form $X({\tt h})=(X({\tt h})^1, X({\tt h})^2, X({\tt h})^3, X({\tt h})^4, X({\tt h})^5)$. Notice that there is no difference between the It\^o- and Stratonovich-form thanks to the specific structure of the diffusion coefficient in our case. 

The support theorem in its classical form is stated for diffusions whose parameters are homogeneous in time. In order to fit into this framework, we replace as before the $5-$dimensional process $X^n$ by a $6-$dimensional process $ (t, X^n _t ) $ which is now a classical time-homogenous diffusion process. This shows that the support theorem applies directly also in the time inhomogeneous case. 
As a consequence, see e.g. Theorem 3.5 of Millet and Sanz-Sol\'e 1994 or Theorem 4 of Ben Arous, Gradinaru and Ledoux 1994, the support of the law $\mathbb{P}^n_{0, x}$ restricted to $ {\cal F}_t$ is the closure 
of the set $ \{ X({\tt h}) : {\tt h} \in {\cal H} \} $ with respect to the uniform norm on $C( [ 0, t], \bbr^5 )  .$  

In order to find lower bounds for (\ref{eq:tobelb}) we have to construct solutions $X({\tt h})$ of (\ref{eq:control1}) which stay in $K_n$ during $ [ 0, t ] . $ On $K_n,$ both processes $X^n$ and $X$ have the same coefficients. Hence, by restricting to $K_n,$ the above control problem (\ref{eq:control1}) is equivalent to  
$$
\left\{\begin{array}{l}
\frac{d}{ds }  X({\tt h})^1_s   \;=\;  \frac{d}{ds} X({\tt h})^5_s \;- F( X({\tt h})^1_s, X({\tt h})^2_s, X({\tt h})^3_s, X({\tt h})^4_s) \\
\frac{d}{ds }  X({\tt h})^2_s \;=\;  \, \al_n(X({\tt h})^1_s)\,(1-X({\tt h})^2_s)  \;-\; \beta_n(X({\tt h})^1_s)\, X({\tt h})^2_s  \,  \\
\frac{d}{ds }  X({\tt h})^3_s\;=\;  \, \al_m(X({\tt h})^1_s)\,(1-X({\tt h})^3_s)  \;-\; \beta_m(X({\tt h})^1_s)\, X({\tt h})^3_s  \,\\
\frac{d}{ds }  X({\tt h})^4_s \;=\; \, \al_h(X({\tt h})^1_s)\,(1-X({\tt h})^4_s)  \;-\; \beta_h(X({\tt h})^1_s)\, X({\tt h})^4_s  \,  \\
\frac{d}{ds }   X({\tt h})^5_s\;=\; (\, S(s)-X({\tt h})^5_s\,)\, \tau  \;+\;  \gamma \sqrt{\tau} \;   \dot {\tt h} (s)  .
\end{array}\right. 
\leqno{\rm (HHcontrolled)}
$$

We construct an explicit solution of ${\rm (HHcontrolled)}$ starting from the fixed initial condition $x = ( v, n , m , h , \zeta ) \in K_{n_0} $ at time $0$ in the following way. First, we choose a path $ \bar v_t = \gamma (t) v ,$ $ t \geq 0 ,$ going from $v$ to $0.$ Here, $ \gamma $ is a smooth function $ \bbr_+ \to [ 0 , 1 ] , \gamma ( 0 ) = 1 , \gamma ( t ) = 0 $ for all $ t \geq 1 .$ Hence for all $ t \geq 1,$ $ \bar v_t \equiv 0 ,$ and $\bar v_0 = v.$ 

Then, solving the equations for $n,$ $m$ and $h$ explicitly, for this fixed choice of $ \bar v_t,$ we obtain
$$ \bar n _t = n e^{ - \int_0^t  a_n ( \bar v_s) ds } + \int_0^t b_n ( \bar v_u) e^{ - \int_u^t a_n ( \bar v_r) dr } du , $$
where $ a_n = \alpha_n + \beta_n , $ $ b_n = \alpha_n .$ We have analogous representations for $\bar m_t $ and $\bar h_t .$  Since $ \bar v_t \equiv 0 $ for all $t \geq 1, $ it follows that 
\begin{equation}\label{eq:expdecay}
| \bar n_{ t + 1 } - n_\infty ( 0 ) | \le C  e^{ - t a_n ( 0 ) } ,
\end{equation}
where the constant depends on $v $ and $n.$ The same convergence result holds for $ \bar m_t $ and $\bar h_t.$ 

Fix $\varepsilon $ such that $ B_\varepsilon ( y_* ) \subset U .$ Then there exists $ t_0$ such that for all $t \geq t_0, $ 
\begin{equation}\label{eq:as}
(\bar n_t, \bar m_t, \bar h_t ) \in B_{ \varepsilon /2 } ( n_\infty ( 0), m_\infty (0), h_\infty (0) ).
\end{equation}
Now we want to choose $ {\tt h}$ such that 
\begin{equation}\label{eq:sol}
\frac{d}{ds }  X({\tt h})^5_s = \frac{d}{ds }  \bar v_s + F ( \bar v_s, \bar n_s, \bar m_s, \bar h_s )  , \quad \mbox{ for all $s \geq  0.$ }
\end{equation}
Equation (\ref{eq:sol}) implies that 
$$ X({\tt h})^5_s =  \zeta + \bar v_s - v + \int_0^s F ( \bar v_u, \bar n_u, \bar m_u, \bar h_u ) du    = : {\tt J^h}_s = {\tt J^h}_s ( v, n, m, h , \zeta ) .$$    
Hence, if we define
$$ \dot {\tt h} (s) := \frac{ \frac{d}{ds }  \bar v_s + F ( \bar v_s, \bar n_s, \bar m_s, \bar h_s )+ ( {\tt J^h}_s - S(s) ) \tau  }{  \gamma \sqrt{\tau} } ,  $$
then $ (\bar v_s, \bar n_s, \bar m_s, \bar h_s, X({\tt h})^5_s)^T $ is indeed a solution of {\rm (HHcontrolled)}, for this specific choice of ${\tt h} .$ 

Fix $t \geq t_0.$ Notice that $\dot{\tt h} $ is well-defined and that $\dot{\tt h} \in L^2 ( [ 0, t ] ) ,$ hence ${\tt h} \in {\cal H} .$ With this choice of ${\tt h},$ the first four lines of ${\rm (HHcontrolled)}$ reduce to the deterministic system (HH)  with input signal $s \to \frac{d}{ds }  \bar v_s + F ( \bar v_s, \bar n_s, \bar m_s, \bar h_s ) .$ Write $\mathbb{Y}$ for the associated deterministic solution starting from $ (v, n,m,h) $ at time $0$ and $\mathbb {X}^x_s = ( \mathbb{Y}_s, {\tt J^h}_s ), s \le t ,$ starting from $x$ at time $0.$ For $n$ sufficiently large, $ \mathbb{X}^x_s \in K_n $ for all $ s \le t .$ By the support theorem, for every $\delta > 0, $ putting $B^\infty_\delta ( \mathbb{X}^x ) = \{ f \in  {\cal C}  : \sup_{ s \le t } | f(s) - \mathbb{X}^x_s | <  \delta \} ,$ we have that 
$$ \mathbb{P}^n_{0, x} (B_\delta^\infty (\mathbb{X}^x) ) > 0 .$$
 
Now, choose $ \delta \le \varepsilon/2  $ and $n$ sufficiently large such that $ B^\infty_\delta ( \mathbb{X}^x )  \subset \{ f \in {\cal C} : T_n (f) > t \} .$ Since $B_\varepsilon ( y_*) \subset U $ and recalling (\ref{eq:as}), we have that $B_\delta ^\infty ( \mathbb{X}^x ) \subset   \{ f \in {\cal C} : f(t) \in B_\varepsilon (y_*) \times \bbr \}    .$ This implies
$$ P_x ( X_t \in B_\varepsilon (y_*) \times \bbr  ) \geq \mathbb{P}^n_{0, x } (B_\delta ^\infty ( \mathbb{X}^x ) ) > 0 ,$$
which finishes our proof. 
\halmos

\section{Proof of (\ref{eq:firststatement}).}\label{sec:proofOU}
For the convenience of the reader we will recall basic concepts from Control Theory as exposed in Sussmann 1973. As above, in order to be able to deal with the $T-$periodic drift coefficient, we work with the space-time process $\bar X_t = ( i_T (t) , X_t) ,t \geq 0 .$ 

Recall that the drift and diffusion coefficients $\bar b(t,x) $ and $\bar \sigma (t,x) $ of $\bar X$ have been introduced in (\ref{eq:b}). We introduce the following family of control vector fields
\begin{equation}\label{eq:D}
G = \{  \bar b + c \bar \sigma  , c \in \bbr \} .
\end{equation}
Here, by definition of $G$, the control parameter $c$ acts on the diffusion part only. There is no control on the drift part. Control vector fields from $G$ correspond to the controls ${\tt h}$ of Section \ref{sec:supp} in case of piecewise constant $\dot {\tt h} .$

Let $G^* $ be the smallest set of vector fields containing $G$ which is closed under Lie brackets. We introduce the mapping $\Delta_{G^*} $ which assigns to every space-time point $ (t,x) \in \bbt \times E_5 $ the linear subspace 
$$ \Delta_{G^*} ( t,x) = Span \{ V^*  ( t, x ) : V^* \in G^* \} .$$
Notice that 
$$ \Delta_{G^*} ( t,x) = Span \{ \bar b(t,x)  ,  L (t, x) : L \in {\cal L}  \} ,$$
where the Lie algebra ${\cal L} $ has been introduced in (\ref{eq:iteration}).  

We say that two points $ ( t , x)$ and $  ( t^* , x^*) $ in $ \bbt \times E_5$ belong to the same orbit of $G$ if and only if there exists a curve $ \gamma $ defined on some interval $ [a,b] $ and a suitable partition $ a = t_0 < t_1 < \ldots < t_r = b $ such that $ \gamma ( a) = (t,x), \gamma (b) = (t^*, x^*)$ and such that on each $ ]t_{i-1}, t_i [  $ there exists a constant $c_i \in \bbr $ with either 
\begin{equation}\label{eq:controlforward}
 \dot \gamma ( t) = \bar b ( \gamma (t)) + c_i \bar \sigma ( \gamma (t))
\mbox{
or   } \dot \gamma ( t) = - \bar b ( \gamma (t)) + c_i \bar \sigma ( \gamma (t)).
\end{equation}

Since the coefficients of $\bar b$ and $\bar \sigma$ are analytic, by Nagano 1966, see also Sussmann 1973, Theorem 8.1 and Section 9, we know that for any $G-$orbit $S \subset \bbt \times E_5 ,$ the following holds.  For all $ (t,x) , (t^*, x^* ) \in S,$ we have that 
\begin{equation}\label{eq:sussmann}
 dim \Delta_{G^*} (t,x) =  dim \Delta_{G^*} (t^*,x^*). 
\end{equation}
In particular, this implies the following. 

Suppose that  $ (0,x) $ and $ (t^*, x^* )$ belong to the same $G-$orbit $S$ such that $ dim \Delta_{G^*} (t^*,x^*) = 6 . $ This is  equivalent to $ dim \Delta_{{\cal L}_N^*} (t^*,x^*) = 5  $ for some $N \geq 1 ,$ where we recall the definition of $\Delta_{{\cal L}_N^*} $ in Definition \ref{def:hoer}. Then we have also for the starting point $ (0, x ) $ the full dimension 
$$ dim \Delta_{G^*} (0,x ) = 6 \quad \mbox{ or equivalently } \quad dim \Delta_{{\cal L}_N^*} (0, x) = 5 \mbox{ for some } N \geq 1. $$

We are now ready to give the proof of (\ref{eq:firststatement}). 

{\bf Proof of (\ref{eq:firststatement}).} 
In the following, we will work with piecewise constant control functions which we call `admissible control'. Our proof relies on the fact that the support of $ \mathbb{P}_{0, x } $ is the closure of all paths $X( {\tt h}) ,$ as defined in (\ref{eq:control1}), where $ \tt h $ is an admissible control.

We prove the following fact. If $X_t \notin {\cal I}_5,$ then $X_{t + s } \notin {\cal I}_5$ for all $s \geq 0$ almost surely. 

Conditioning on $X_t = x ,$ we can assume without loss of generality that $t = 0 $ and $  X_0 = x \notin {\cal I}_5.$ Thus, all deterministic control paths $X({\tt h})$ issued from $x$ and using an admissible control are such that the curve $ ( i_T (t), X({\tt h})_t)$ belongs to an orbit of $G$ on which the dimension of $\Delta_{G^*} $ is strictly less than $6.$ This implies that for any fixed $N,$ $dim \Delta_{{\cal L}_N^*} <  5 $ on the whole orbit. Since this holds for any $N,$ this implies that $ X_s  \notin  {\cal I}_5 $ for all $ s \geq 0 ,$ $ P_x-$almost surely. This concludes our proof. \halmos

\section{Proof of Theorem \ref{theo:ergodicOU} Part c). }\label{sec:arnold-kliemann}
We still work with the space-time process $\bar X_t = ( i_T (t) , X_t) ,t \geq 0,$ taking values in $\bbt \times E_5 .$ The process $ \bar X_t $ has the transition operator $\bar P_t $ given by 
$$ \bar P_t ((s,x) , \cdot ) = \delta_{ i_T ( t+s) } \otimes P_{ s, s+ t } ( x, \cdot ) .$$ 
We shall denote invariant probability measures of $\bar X_t $ by $ \bar \mu .$  

In order to prove Part c) of Theorem \ref{theo:ergodicOU}, we use control sets as in Arnold and Kliemann 1987, to characterize the support of extremal invariant probability measures $\bar \mu $ of $\bar X_t .$ For that sake, for any $ (s, x ) \in \bbt \times E_5  $ and any $ t > 0 , $ we put 
\begin{multline*}
 {\cal O}^+ ( t, (s,x) ) = \{ \gamma (t) :  \mbox{ there exists an admissible control $ {\tt h}$ such that  }\\
 \gamma ( s) = x + \int_0^s [ \bar b ( u, \gamma (u)) + \bar \sigma ( \gamma (u))\dot {\tt h} ( u)] du \; \mbox{ for all } \; s \le t \} .
\end{multline*} 
Notice that in the above definition we are moving through the orbit forward in time. In other words, $ {\cal O}^+ ( t, (s,x) )$ is the set of all points reachable from $ (s,x) $ forward in time during a time period of length $t.$ We will also note
$$ {\cal O}^+ ( s, x) = \bigcup_{ t > 0 } {\cal O}^+ ( t, (s,x) ) ,$$
the set of all points reachable forward in time, starting from $(s,x).$ Then a set $ F \subset \bbt \times E_5 $ is called an {\it invariant control} set if 
$$ \overline{ {\cal O}^+ ( (s,x) ) } = \bar F \mbox{ for all } (s,x) \in F .$$
Notice that invariant control sets are necessarily disjoint. By Proposition 1.1 of Arnold and Kliemann 1987, to any extremal invariant probability measure $\bar \mu $ is associated a unique invariant control set $F$ such that $ supp \; \bar \mu = \bar F .$ 

In the following we start by describing the relationship between extremal invariant probability measures $ \bar \mu$ of the space-time process $ \bar X_t$ and extremal invariant probability measures $\mu$ of the skeleton chain $ (X_{kT } )_{ k \geq 0 }.$ Then we prove that there are only finitely many invariant control sets.

\begin{prop}\label{prop:prop4}
The following assertions are equivalent. 
\begin{enumerate}
\item 
$\mu $ is an invariant probability measure of the skeleton chain $ (X_{kT})_{k \geq 0 } .$
\item
For any $s \in ] 0, T [, $ the measure $ \mu_s := \mu P_{ 0, s } $ is an invariant probability measure of  $(X_{k T + s })_{ k \geq 0 } ,$
and $ \mu_s P_{s, s+ t } = \mu_{i_T  (s+t) } .$ 
\item
The measure
\begin{equation}\label{eq:barmu}
 \bar \mu := \frac1T \int_0^T ds ( \delta_s \otimes \mu_s ) 
\end{equation}
is an invariant probability measure of $ \bar X .$ 
\end{enumerate}
\end{prop}

{\bf Proof of Proposition \ref{prop:prop4}.}
It is straightforward to show the equivalence of the first two points 1. and 2.  Formula (\ref{eq:barmu}) then shows how to build invariant measures for the process $ \bar X$ starting from invariant measures of the skeleton chain. 

We have to show that any invariant measure $ \bar \mu$ for $ \bar X$ can be written in a form (\ref{eq:barmu}). The first marginal $ \bar \mu_{ | t } $ is necessarily the uniform measure on the torus $\bbt. $ Then by Lebesgue disintegration we have $ \bar \mu = \frac1T \int_0^T ds ( \delta_s \otimes K(s, \cdot ) ),$ where $K(s, dx )$ is a regular version of the conditional distribution of the second component of $\bar\mu$ given the first component. By invariance,  $\bar \mu \bar P_h = \bar \mu$ for all $h>0$. Take $\bar \mu$ as starting law for $\bar X$. If value $s$ has been selected for the first component, then $\tilde\mu:= K(s,\cdot)$ acts as starting law for the second component. The construction of $\bar X$ and the explicit form of the transition $\bar P_h $ yield $K(s,\cdot) P_{s,s+h}=\call_{\tilde\mu}(X_h)$. Under the same starting condition, the first component of $\bar X _h$ is $i_T(s+h)$, thus  the second component of $\bar X _h$ has law $K(i_T(s+h),\cdot)= \call_{\tilde\mu}(X_h)$.

Note that for every law $\wt\mu$ on $E_5$ and every $f\in \calc_b(E_5)$,  $\,h\to E_{\tilde\mu}(f(X_h))\,$ is continuous, by continuity of the sample paths of $X$. In particular, taking $ \wt\mu = K(s, \cdot ) $ as above and $ h = T-s + t$, $0\le t \le T, $ we obtain that $ t \to \int_{E_5} f(x)  K( t, dx ) $ is continuous. This implies that we can take a version of $K(\cdot,\cdot)$ such that $t\to K(t,\cdot)$ is continuous. 

This avoids problems related to $\la(ds)-$null sets in the conditional expectations. Thus we have proved that $K(s,\cdot)P_{s,s+h} =K(i_T(s+h),\cdot)$ for all $s\in \bbt$ and all $h>0$.

$\bbt$ being the torus, invariance $\bar \mu \bar P_T = \bar \mu$ now gives $K(s,\cdot)=K(s,\cdot) P_{s,s+T}$  for all $s$ in $\bbt$. Thus $\mu_s:=K(s,\cdot)$ is an invariant measure for $(X_{kT+s})_{k\ge 0}$, and with $\mu:=K(0,\cdot)$ we have 1. and 2.  \halmos\\

The following proposition follows easily from the above considerations. 

\begin{prop}\label{prop:prop5}
$\bar \mu = \frac1T \int_0^T ds ( \delta_s \otimes \mu_s ) $ is an extremal invariant measure of $\bar X $ if and only if  any $\mu_s $ is an extremal invariant measure of  $(X_{kT + s })_{ k \geq 0 } .$  
\end{prop}

{\bf  Proof of Proposition \ref{prop:prop5}.} Suppose that $ \mu_s = \mu P_{ 0, s } $ is extremal for $(X_{kT + s })_{ k \geq 0 } .$ We have to show that  $ \bar \mu = \frac1T  \int_0^T ds ( \delta_s \otimes \mu_s ) $ is extremal for $\bar X.$  Suppose that there exists $ \alpha \in ] 0, 1 [ $ such that 
$$ \bar \mu = \alpha \bar \mu_1 + ( 1 - \alpha ) \bar \mu_2 $$
with $ \bar \mu_1 \neq \bar \mu_2 $ invariant measures of  $ \bar X.$  Lebesgue disintegration yields 
$$ \bar \mu_i = \frac1T  \int_0^T ds ( \delta_s \otimes \mu^i _s ) , i = 1 , 2 ,$$
and therefore
$$ \bar \mu = \frac1T \int_0^T ds ( \delta_s \otimes ( \alpha \mu^1_s + ( 1 - \alpha ) \mu^2_s ) = 
\frac1T \int_0^T ds ( \delta_s \otimes \mu_s ) .$$
This implies that for all $ s \in \bbt,$ 
$ \mu_s = \alpha \mu^1_s + ( 1 - \alpha ) \mu^2_s ,$ where $ \mu^1_s $ and $ \mu^2_s $ are invariant measures of the skeleton chain $ (X_{kT + s })_{ k \geq 0 } .$ Since $\mu_s $ is extremal,  it follows from this that 
$ \mu_s = \mu^1_s = \mu^2 _s , $ for all $ s \in \bbt, $ which implies that $ \bar \mu = \bar \mu_1 = \bar \mu_2 , $ which is a contradiction. On the other hand it is straightforward to show that $\bar \mu$ extremal implies that any $\mu_s$ is an extremal invariant measure.

\halmos

We are now able to give the proof of Theorem \ref{theo:ergodicOU} c). We have already shown that the skeleton chain possesses invariant probability measures. We now show that the skeleton chain possesses only a finite number of extremal invariant probability measures. Let $\mu$ be such an extremal measure and let $ \bar \mu $ be the associated extremal invariant measure of $ \bar X.$ Then there exists an invariant control set $ F $ such $supp \; \bar \mu= \bar F .$ Fix a starting point $ (0, x ) \in F  $ and consider the process issued from $ (0, x) .$ Then the skeleton chain $ ( X_{ k T } )_{ k \geq 0 } $ starting from $x$ at time $0$ induces a subset $ \{ ( 0, X_{ k T } ) , k \geq 0 \} \subset  \{ \bar X_t : t \geq 0 \} ;$ it is this subset which is in the center of our interest. 

With the notation of Proposition \ref{prop:prop4}, let $A_k$ denote the set of all paths which visit ${B_{\varepsilon_k } (x_k)} $ infinitely often, $ A_k' $ the set of all paths which visit ${B_{\varepsilon_k } (y_k)} $ infinitely often.
Since $K$ is visited i.o. almost surely, there exists an index $ k \in \{ 1 , \ldots , N \} , $  such that $P_x ( A_k) > 0 .$ Nummelin splitting then shows that $P_x ( A_k') > 0 .$ This means that $ B_{ \varepsilon_k } ( y_k ) $ belongs entirely to the support of $ \sum_{ k \geq 1 } e^{- k } P_{ 0, k T } ( x, \cdot ) ;$ i.e., 
\begin{equation}\label{eq:schoen}
 B_{ \varepsilon_k } ( y_k ) \subset supp \; \sum_{ k \geq 1 } e^{- k } P_{ 0, k T } ( x, \cdot )  .
\end{equation}
But by the support theorem, 
$$ supp \; \sum_{ k \geq 1 } e^{- k } P_{ 0, k T } ( x, \cdot ) = \overline{\bigcup_{ k \geq 1 } \Pi_2\left( {\cal O}^+ ( kT , (0,x))\right) } ,$$
where $ \Pi_2$ denotes the projection on the space variable. Thus, using (\ref{eq:schoen}) and the fact that $F$ is an invariant control set,
$$ \{ 0 \} \times B_{ \varepsilon_k } ( y_k ) \subset \overline{\bigcup_{ k \geq 1 }  {\cal O}^+ ( kT , (0,x)) } 
\subset \overline{ {\cal O}^+ ( (0,x))} = \bar F .$$
Hence any invariant control set $ F$ which is the support of an extremal invariant measure $\bar \mu $ is such that its closure contains (at least) one of the finitely many balls $ \{ 0 \} \times B_{\varepsilon_k} ( y_k ) .$ Since invariant control sets are pairwise disjoint, there are no two control sets that can contain the same ball $ \{ 0 \} \times B_{\varepsilon_k} ( y_k ) $ at the same time. Thus there exist only finitely many such invariant control sets, that is, only finitely many extremal invariant probability measures $ \bar \mu $ of $\bar X,$ hence by Lebesgue disintegration, only finitely many extremal invariant probability measures $ \mu $ of the skeleton chain. This concludes our proof. \halmos

\section*{Acknowledgments}
We thank Michel Bena\"im very warmly for stimulating discussions on control arguments. We also thank two anonymous referees for helpful comments and suggestions. 

\bibliography{biblio}

\nocite{arnold-kliemann}
\nocite{aihara}
\nocite{benarous}
\nocite{berglund1}
\nocite{berglund2}
\nocite{hoepfnerbrodda}
\nocite{desroches}
\nocite{destexhe}
\nocite{guckenheimer}
\nocite{endler}
\nocite{hodgkin-huxley}
\nocite{hoepfner2007}
\nocite{ho-ko}
\nocite{h-l-t}
\nocite{izhi}
\nocite{ikeda-wa}
\nocite{K-S}
\nocite{kawazu-watanabe}
\nocite{kliemann}
\nocite{colonius}
\nocite{Kunita}
\nocite{Kunita2}
\nocite{Kusuoka-Stroock}
\nocite{M-T}
\nocite{Millet}
\nocite{Morris}
\nocite{david}
\nocite{nummelin}
\nocite{nummelin2}
\nocite{pankratova}
\nocite{rinzel-miller}
\nocite{rubin-wechselberger}
\nocite{stroock-varadhan}
\nocite{sussmann}
\nocite{wang}

\end{document}